\documentclass[a4paper, 10pt, conference]{ieeeconf}      

\IEEEoverridecommandlockouts                              

\usepackage{graphicx}      
\usepackage{amsmath,amssymb,amsopn}
\usepackage{dsfont}
\usepackage{subcaption} 

\usepackage{mathtools}
\mathtoolsset{showonlyrefs}

\newtheorem{Remark}{Remark}

\newenvironment{remark}[1][]{\begin{Remark}}{\hfill$\square$\end{Remark}}

\title{\LARGE \bf Numerical Investigation of Traffic State Reconstruction and Control Using Connected Automated Vehicles}

\author{Mladen \v{C}i\v{c}i{\'c}$^{1}$, Matthieu Barreau$^{1}$ and Karl Henrik Johansson${}^1$
\thanks{${}^1$ Division of Decision and Control Systems, KTH Royal Institute of Technology Stockholm, Sweden (e-mail: cicic,barreau,kallej@kth.se).}%
\thanks{The research leading to these results has received funding from the KAUST Office of Sponsored Research under Award No. OSR-2019-CRG8-4033, the European Union's Horizon 2020 research and innovation program under the Marie Sk{\l}odowska-Curie grant agreement No 674875, VINNOVA within the FFI program under contract 2014-06200, the Swedish Research Council, the Swedish Foundation for Strategic Research and Knut and Alice Wallenberg Foundation. The authors are affiliated with the Wallenberg AI, Autonomous Systems and Software Program (WASP).}
}

\begin{document}

\maketitle
\thispagestyle{empty}
\pagestyle{empty}

\begin{abstract}
In this paper we present a numerical study on control and observation of traffic flow using Lagrangian measurements and actuators. We investigate the effect of some basic control and observation schemes using probe and actuated vehicles within the flow. The aim is to show the effect of the state reconstruction on the efficiency of the control, compared to the case using full information about the traffic. The effectiveness of the proposed state reconstruction and control algorithms is demonstrated in simulations. 
They show that control using the reconstructed state approaches the full-information control when the gap between the connected vehicles is not too large, reducing the delay by more than $60\%$ when the gap between the sensor vehicles is $1.25$~km on average, compared to a delay reduction of almost $80\%$ in the full-information control case. 
Moreover, we propose a simple scheme for selecting which vehicles to use as sensors, in order to reduce the communication burden. Numerical simulations demonstrate that 
with this triggering mechanism, the delay is reduced by around $65\%$, compared to a reduction of $72\%$ if all connected vehicles are communicating at all times.
\end{abstract}

\section{Introduction}

Traffic flow control has attracted a lot of attention in the past decades. Indeed, the current traffic infrastructure cannot handle the increasing number of cars and trucks, leading to congestion, which in turn wastes productivity, increases the emissions and jeopardizes the safety of road users. 
Since it is not always possible, nor desirable, to add lanes or create auxiliary roads, it is imperative to have control over the traffic flow.
Good awareness of the traffic situation is crucial for implementation of any form of traffic flow control.

In urban environments, the complex structure of the road network makes modelling the system and controlling its traffic flow difficult.
Even in a freeway context, where modelling is easier, there are still difficulties due to the high dimensionallity and nonlinearity of the discretized model of the road, among other reasons.
Nevertheless, there exist many control schemes based e.g. on model predictive control \cite{ferrara2018state}. 

To avoid discretization, one may use the infinite-dimensional nonlinear Partial Differential Equation (PDE) model \cite{lighthill1955kinematic,richards1956shock}. In the case of a Greenshields' fundamental diagram \cite{lighthill1955kinematic}, the obtained PDE is hyperbolic and control is usually applied at the boundaries of the PDE, i.e. at entrance and exit, using ramp metering for instance \cite{ferrara2018state}. The control strategies can be of various types: explicit state feedback \cite{karafyllis2017analysis}, PI control \cite{zhang2019pi} or backstepping \cite{yu2018traffic}. Nevertheless, boundary control is difficult to derive and may lead to unrealistic infinite-dimensional control laws. 
Recently, there has been a push for in-domain, or Lagrangian traffic control.
In \cite{PIACENTINI201813}, the authors use a moving bottleneck to help dissipate congestion faster. A similar idea is used in \cite{cicic2018traffic,cicic2019multiclass}.

Almost all the traffic control laws developed previously require the knowledge of the full-state. In other words, either there are sensors everywhere along the road or the state is estimated using an observer. There is not much work on observers for traffic flow (see the survey in \cite{seo2017traffic}). In the recent article \cite{YU2019183} for instance, a backstepping observer is derived for a congested road only. 
However, from the point of view of traffic control, the most interesting case is when we have a mixture of congestion and free flow.
In \cite{qi2018boundary}, an estimate of the flow in this case is considered but it relies on a merger of the free flow and congested flow estimates.
More recently, in \cite{seo2015probe,herrera2010incorporation,amin2008mobile,delle2019traffic}, observers based on measurements from probe vehicles are proposed. Using vehicles as sensors makes sense since the vehicles today are well-equipped and, within a decade, they may be able to communicate with the infrastructure.

The main contribution of this paper is continuing the work started in \cite{cicic2018traffic}, using Connected Automated Vehicles (CAVs) as moving bottlenecks to dissipate stop-and-go waves. The main novelty is that instead of assuming we have full information of the traffic conditions everywhere on the road, we base the control on the reconstructed traffic state which is obtained using probe vehicles as sensors. The proposed algorithm for reconstruction is designed so that the communication burden is reduced, using a triggering mechanism to activate probe vehicles where additional sensing is required.

To this end, in Section~\ref{sec:model}, the models used in the article are stated. Next, in Section~\ref{sec:control} a control law for stop-and-go wave dissipation using controlled moving bottlenecks is presented.
Then, in Section~\ref{sec:observe} we describe a simple algorithm, used to reconstruct the traffic density on the road, based on local probe vehicle measurements.
Finally, in Section~\ref{sec:simulations} we test the state reconstruction and control schemes and discuss the simulation results, and in Section~\ref{sec:conclude}, we conclude.

\section{Model for Observation and Control}
\label{sec:model}

Throughout this paper, we are interested in modeling a portion of a freeway with two lanes. We adopt the discretization of the Lighthill-Whitham-Richards model \cite{lighthill1955kinematic},
\begin{equation}
\label{eq:LWR}
	\partial_\tau \rho (\tau,x) + \partial_x \left( v \rho \right)(\tau,x) = 0
\end{equation}
where $\rho(\tau,x) \in [0, P]$ is the average density on the road at point $x$ and time $\tau$ over the two lanes. The average speed of the vehicles at some point $v(\tau,x)$ is directly determined by the traffic density.
Alternatively, it can be easier to express the evolution of the traffic density using the traffic flow $q(\tau,x) = \rho(\tau,x) v(\tau,x)$, representing the number of vehicles passing through point $x$ at time $\tau$, and can usually be directly measured using fixed sensors such as induction loops.
There are many different ways to model the dependence of the flow on traffic density, i.e. the fundamental diagram.
In this work, we will be using the triangular, Newell-Daganzo fundamental diagram, 
\begin{equation}
	q(\rho) = \min\left( V \rho, \ Q, \ W(P-\rho) \right),
\end{equation} 
where $P$ is the maximal density on the road, $V$ is the free flow speed, and $W$ is the congestion wave speed.

In what follows, we first present the numerical scheme that will be used to simulate \eqref{eq:LWR}, and then briefly describe how the model can be extended to capture the behaviour of stop-and-go waves and moving bottlenecks.

\subsection{Extended CTM}

A standard way of numerically solving and simulating \eqref{eq:LWR} is by applying Godunov discretization, which yields the well-known Cell Transmission Model (CTM) \cite{DAGANZO1994269}.
Here, we use the extended version of the CTM similar to the one in \cite{cicic2019multiclass}, where the update of traffic density in cell $i$, $\rho_i(t)$, is
\begin{align}
    \label{eq:CTM}
    \rho_i(t+1) &= \rho_i(t) + \frac{T}{L}\left(q_{i-1}(t) - q_{i}(t)\right),\\
    q_i(t) &= \min\left(U_i(t)\rho_i(t), Q_i(t), W\!\left(P \!- \rho_{i+1}(t)\right)  \right),
\end{align}
for each of $i=1, \ldots, N$ cells, and $q_0(t)$ is the inflow to the road segment that is defined externally.
The cell length $L$ and time step $T$ are linked by $L=VT$.
This model differs from the standard one in that we model the capacity drop phenomenon \cite{ferrara2018state} through letting
\begin{equation}
    Q_i(t) = \min\left(V \sigma, W\left( P-(1-\alpha)\sigma - \alpha\rho_i(t)\right)\right),
\end{equation}
where $\sigma \in [0, P]$ is the critical density and $\alpha \in [0, 1]$ and by allowing the free flow speed of vehicles in cell $i$ to vary in time. For $\alpha = 1$, we recover the previous case.

Given a reference traffic density profile $\rho_i^k(t)$ that should at times $t$ be satisfied in cells $i_-^k(t)$ to $i_+^k(t)$, we may set
\begin{equation}
\label{eq:Uhead}
    U_{i}(t) = V\min\left\{1, \frac{\rho_{i+1}^k(t+1)}{\rho_i(t)}\right\},
\end{equation}
for $i = i_+^k(t)-1$, and then recursively have
\begin{equation}
\label{eq:Ubody}
\resizebox{0.91\columnwidth}{!}{$
    \!\!\!U_i(t) \!=\! V\!\min\!\left\{\!1,\max\!\left\{0,\!\!\frac{\rho_i^k(t+1)-\frac{V-U_{i+1}(t)}{V}\rho_{i+1}(t)}{\rho_i(t)}\!\right\}\!\!\right\},
    $}
\end{equation}
for $i = i_+^k(t)-2$ down to $i_-^k(t)-1$, in order to make the actual traffic density converge to the reference.
Where not otherwise noted, we use the default free flow speed, ${U_i(t) = V}$.
Since the evolution of the reference density profile in time is known for cells around a stop-and-go waves and moving bottlenecks, this enables proper handling of these phenomena, which would otherwise suffer from diffusion that is inherent in the classical CTM.
Reference density profiles will be described in the following subsections.

\subsection{Stop-and-go waves}

In this work, we focus on stop-and-go waves that originate from downstream of the road segment of interest.
Once the $k$--th stop-and-go wave, with density $\rho_c^k>\sigma$, stretching from cell $i_t^k(t)$ to cell $i_h^k(t)$ and with downstream end at $z^k(t)$, is fully within the road segment, traffic density in its zone should follow \cite{cicic2019multiclass}:
\begin{equation}
    \!\rho_i^k(t) \!= \!\begin{cases}
    \rho_c^k,&\!\!\!\! i=i_t^k(t), \ldots, i_h^k(t),\!\!\!\!\!\\
    \rho_d^k\!+\!\left(\rho_c^k \!\!-\!\! \rho_d^k\right)\frac{ z^k(t) - i_h^k(t)L }{L},&\!\!\!\! i=i_h^k(t),\\
    \rho_d^k,&\!\!\!\! i=i_h^k(t)+1,
    \end{cases}
\end{equation}
until $i_h^k(t)=i_t^k(t)$ and the stop-and-go wave is dissipated.
Here, the traffic density discharged from the stop-and-go wave is
\begin{equation}
    \rho_d^k = \frac{W}{V}\left(P-(1-\alpha)\sigma - \alpha \rho_c^k\right),
\end{equation}
and the downstream end of the wave moves according to
\begin{equation}
        z^k(t+1) = z^k(t)+\lambda_d T = z^k(t) - V\frac{(1-\alpha)\sigma}{P-(1-\alpha)\sigma}T.
\end{equation}
This traffic density reference profile holds from $i_-^k(t) = i_t^k(t)$ to $i_+^k(t) = i_h^k(t)+1$.



\subsection{Behavior of CAVs}


Apart from modelling the dynamics of aggregate traffic, we also need to describe the movement of specific CAVs, as well as capture how changing their behaviour influences the rest of the traffic. 
The reason for particular treatment of CAVs is that we assume that they are able to communicate with the infrastructure, sending information about the traffic in their vicinity, and potentially receiving control actions to apply, acting as both pointwise Lagrangian sensors and actuators.
We differentiate between three types of CAVs:
\begin{enumerate}
    \item Inactive CAVs -- not acting as sensors or actuators,
    \item Probe CAVs -- acting only as sensors, and
    \item Actuator CAVs -- acting as both sensors and actuators.
\end{enumerate}

\begin{remark}
Having both probe and actuators CAVs is a more reasonable assumption than considering only actuator CAVs. Indeed, there already exist cars capable of sensing the traffic situation and in the near future, there might be more probe vehicles than fully automated ones. 
The third category reflects the fact that we are not using all CAVs at any times and that we can switch off some of them to decrease the communication burden.
\end{remark}

We denote the position of the $m$--th CAV as $y^m(t)$, and update it according to
\begin{equation}
    y^m(t+1) = y^m(t) + v_y^m(t)T,
\end{equation}
until the CAV reaches the end of the road segment and leaves it.
Initially, we index the vehicles on the road at $t=0$ so that ${y^m(t)\ge y^{m+1}(t)}$, i.e. the lowest $m$ corresponds to the downstream-most CAV, and each of the newly arrived CAV will get a higher index.
This ordering may change in the event of CAVs overtaking each other.
However, the ordering will always be preserved in case of actuator vehicles, ${y^m(t)> y^{m'}(t)}$ if ${m'>m}$, ${m \in \mathcal{Y}_a}$, ${m' \in \mathcal{Y}_a}$, since actuator vehicles will not overtake each other.
The speed of vehicle $m$ is given by
\begin{equation}
    v_y^m(t) = \min \left\{ v_{i^m_y(t)}(t), u^m(t)\right\}
\end{equation}
where $v_{i^m_y(t)}(t)$ is the density-dependent speed of the traffic in cell ${i^m_y(t) = \left\lfloor y^m(t)/L\right\rfloor}$ where the $m$--th CAV is, ${v_i(t) = q_i(t)/\rho_i(t)}$,
and $u^m(t) \in \left[u_{\min},V\right]$ is the control input.
If the $m$--th CAV is inactive, a probe vehicle, or an actuator vehicle not being actively controlled at time $t$, we set ${u^m(t) = V}$ and the vehicle moves with the flow of traffic.

If we externally impose speed ${u^m(t)<v_{i^m_y(t)}(t)}$, i.e. force the $m$--th CAV to move slower than the rest of the traffic, it will begin to affect the traffic flow at its position by acting as a moving bottleneck.
This phenomenon may be modelled by using \eqref{eq:Uhead}, \eqref{eq:Ubody} with density profile reference \cite{cicic2019multiclass}:
\begin{equation}
    \!\rho_i^m(t) \!= \!\begin{cases}
    \rho_b^m(t),&\!\!\!\! i=i_y^m(t)-1\\
    \sigma_b\!+\!\left(\rho_b^m(t) \!\!-\!\! \sigma_b\right)\frac{ y^m(t) - i_y^m(t)L }{L},&\!\!\!\! i=i_y^m(t),\\
    \sigma-\sigma_b,&\!\!\!\! i=i_y^m(t)+1,
    \end{cases}
\end{equation}
and $i_-^m(t) = i_y^m(t)-1$, $i_+^m(t) = i_y^m(t)+1$, and where the density of the congestion in the wake of the moving bottleneck is
\begin{equation}
    \rho_b^m(t) = \frac{WP-\left(V-u^m(t)\right)\left(\sigma-\sigma_b\right)}{u^m(t)+W}.
\end{equation}
For the purpose of control, this congestion is also counted as a stop-and-go wave, with $\rho_d^k = \rho_b^m(t)$, $\lambda_d = u^m(t)$ and $z^k(t) = y^m(t)$.

\section{CAV-based control}
\label{sec:control}

In nature, the density of vehicles on the road is bounded and then traffic flow systems are stable in the classical sense. The aim of the control is therefore to improve the performances. In \cite{ferrara2018state}, many cost functions are introduced to optimize ecological or economical indices. In this work, we focus on the Total Time Spent ($TTS$) defined as the time it takes for cars to exit the road, expressed by
\begin{equation}
    TTS = \sum_{t=1}^{t_{max}}T\left(n_0(t) + \sum_{i=1}^{i_{max}} \rho_i(t)L\right),
\end{equation}
where $n_0(t)$ is the number of vehicles queuing to enter the road at its upstream end.

For the control, we use a subset $\mathcal{Y}_a$ of activated CAVs as actuators for the traffic flow. They are acting as controlled moving bottlenecks, as in \cite{cicic2018traffic,cicic2019multiclass}.
Due to capacity drop, the discharging flow from a stop-and-go wave will be lower than the road capacity, causing an increase $TTS$ for the vehicles on the road.
By reducing the speed of the actuator CAV, we restrict the inflow to the stop-and-go wave to a value smaller than its discharging flow.
As a result, the length of the stop-and-go wave will decrease over time until it is fully dissipated, at which point the capacity of the road is returned to its maximum value.
Each actuator vehicle focuses on dissipating one stop-and-go wave as soon as possible, if there are any stop-and-go waves downstream of them.
If we predict the actuator vehicle $m$ will fail to dissipate its focus stop-and-go wave $k_m$, due to limitations on its minimum speed, we force the next actuator vehicle upstream of the $m$--th CAV to focus on the same stop-and-go wave, $k_{m'} = k_m$.
Otherwise, each actuator vehicle will focus on the first stop-and-go wave downstream of its position, $i_h^{k_m}>i_y^m$.

Denoting by $n^m_{yz}(t)$ the number of vehicles between the $m$--th CAV and the downstream end of the stop-and-go wave $k_m$, we will have
\begin{equation}
    \dot{n}^m_{yz}(t) = \left(V - u^m(t)\right)\left(\sigma-\sigma_b\right) \left(V-\lambda_d\right) \rho_d^{k_m},
\end{equation}
whereas the distance between $y^m(t)$ and $z^{k_m}(t)$ will follow
\begin{equation}
    \dot{d}^m(t) = \dot{z}^{k_m}(t)-\dot{y}^m(t) =  \lambda_d - u^m(t).
\end{equation}
We minimize the adverse effects that the stop-and-go wave has on the traffic flow by ensuring that we have ${n^m_{yz}(\theta) = 0}$ and ${d^m(\theta) = 0}$ with minimum $\theta$, since in that case the road capacity is returned to its maximum value as fast as possible, without excessively delaying the traffic.

This is achieved by setting the speed of vehicle $m$ to
\begin{equation}
\label{eq:um}
    u^m(t) = \frac{V \left(\rho_d^{k_m}-\sigma + \sigma_b\right)- \lambda_d \left( \bar{\rho}_{i_y^m, i_h^{k_m}}(t)-\rho_d^{k_m} \right)}{\bar{\rho}_{i_y^m, i_h^{k_m}}(t) - \sigma + \sigma_b},
\end{equation}
where $\bar{\rho}_{i_y^m, i_h^{k_m}}(t)$ is the average traffic density between the actuator vehicle cell $i_y^m(t)$ and the downstream end of the stop-and-go wave cell $i_h^{k_m}(t)$.
If the speed thus calculated is lower than the allowed minimum,  $u^m(t) < u_{\min}$, we instead apply $u^m(t) = u_{\min}$ and conclude that the $m$--th CAV will not succeed in dissipating the ${k_m}$--th stop-and-go wave.

Note that the initial average densities $\bar{\rho}_{i_y^m, i_h^{k_m}}(t)$ will be calculated using the reconstructed traffic density profile,
\begin{equation}
    \bar{\rho}_{i_y^m, i_h^{k_m}}(t) = \frac{1}{i_h^{k_m}-i_y^m}\sum\limits_{i=i_y^m}^{i_h^{k_m}} \hat{\rho}_(t)
\end{equation}
based on the information that we have available in the particular case.
Therefore, the different cases of control that are considered in this work are only distinguished based on what information we use for traffic state reconstruction.

\section{Traffic State Reconstruction}
\label{sec:observe}

As mentioned in the previous section, in order to be able to improve the traffic flow, we first need to sense and at least approximately reconstruct the traffic density profile along the road.
Here we propose a simple traffic state reconstruction scheme, and briefly discuss how we select which CAVs are used as probe vehicles.

\subsection{Methodology}

We assume that we only have access to the information about the traffic that is communicated by probe and actuator vehicles on the road, from set $\mathcal{Y}_s(t)$.
Namely, we assume that these CAVs can measure local traffic densities in cells ${i \in \mathcal{I}_s(t)}$ adjacent to the cell they are in,
\begin{equation}
\label{eq:adjacentset}
    \mathcal{I}_s(t) = \left\{i: \left|i-i_y^m(t)\right|\le 1, m \in \mathcal{Y}_s(t)\right\},
\end{equation}
assuming that the cell length $L$ is chosen so that the sensors on CAVs can indeed cover this range.
This set will typically change every time step, since the CAVs will move along the road, leave the road segment at its downstream end, and new ones will arrive at its upstream end.
Using these measurements, we can attempt to approximately reconstruct the traffic density,
\begin{equation}
    \hat{\rho}_i(t) = \begin{cases}
    \rho_i(t),& i \in \mathcal{I}_s(t)\\
    \hat{\rho}_i(t\!-\!1) + \frac{T}{L}\left(\hat{q}_{i\!-\!1}(t\!-\!1) - \hat{q}_{i}(t\!-\!1)\right),& i \notin \mathcal{I}_s(t) 
    \end{cases}
\end{equation}
where $\hat{q}_i(t)$ is defined the same way as $q_i(t)$, but using $\hat{\rho}_i(t)$ instead of $\rho_i(t)$, and $\hat{U}_i(t)$ that is calculated for $\hat{\rho}_i(t)$.

Since here we assume that the traffic flow model is known, there are only three sources of uncertainty in the traffic density estimate: the initial condition $\rho_i(0)$ (which disappear in time \cite{delle2019traffic}), the inflow $q_0(t)$, and the conditions at the downstream end of the road segment, i.e. stop-and-go waves arriving from downstream.
We assume that at least the average inflow $\bar{q}_0$ is known, which in practice could be learned from historical data.
Then, we may use this value as the estimated inflow, ${\hat{q}(t) = \bar{q}_0}$, if no other information is available, as well as for initializing the traffic density estimates, ${\hat{\rho}_i(0) = \frac{\bar{q}_0}{V}}$.
However, unless they can be measured in some other way, changes in the traffic conditions downstream of the road segment will only be detected once a probe vehicle reaches their spillback.

\begin{figure}
	\centering
	\includegraphics[width=8cm]{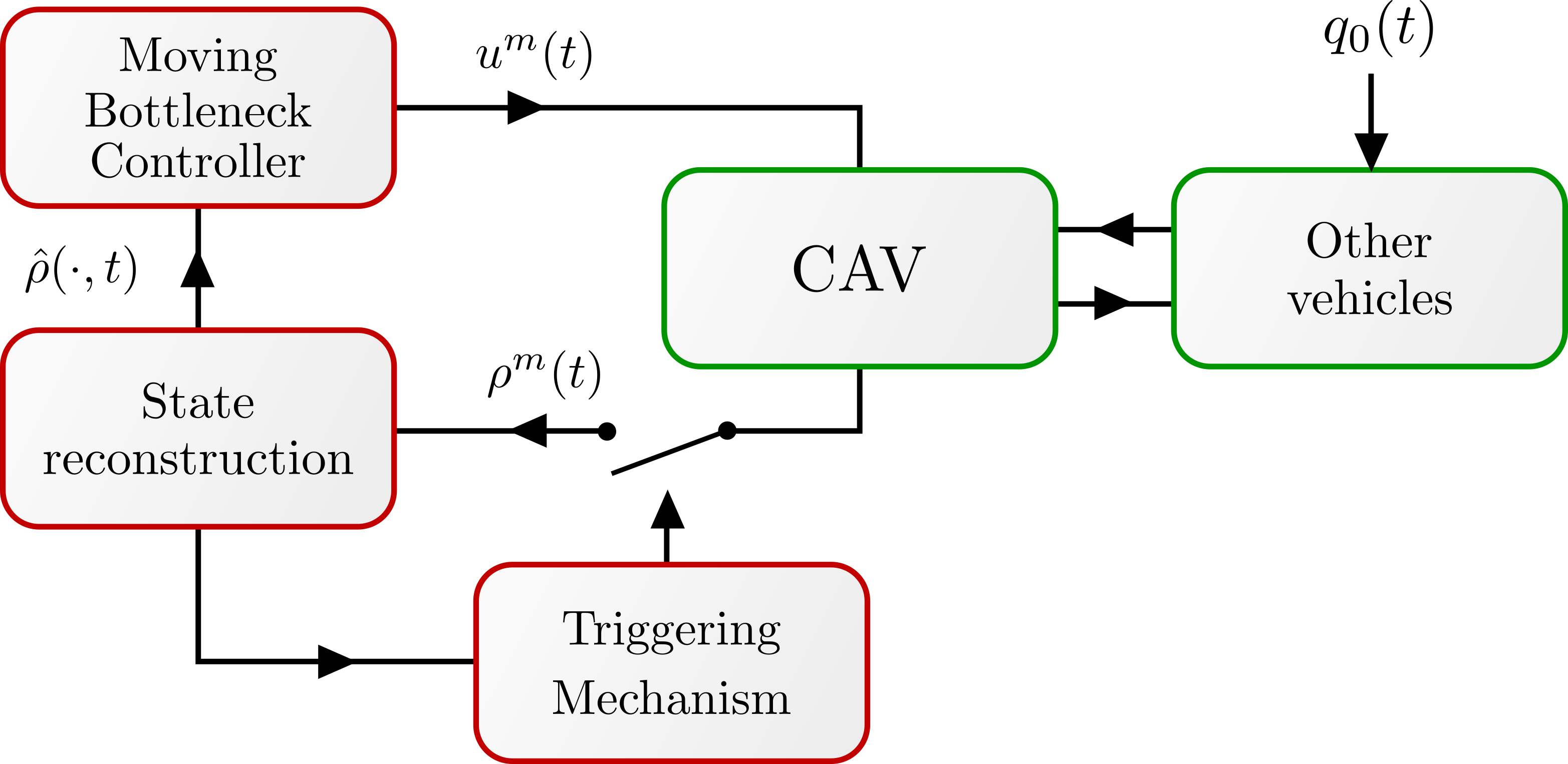}
	\caption{Block diagram of the control scheme used here. The triggering mechanism reduces the communication burden by only activating the probe CAVs when their sensing is needed, based on the reconstructed density.}
	\label{fig:blockDiagram}
\end{figure}

Note that in the proposed simple reconstruction algorithm it is required that all probe vehicles communicate their measurements at each time step, potentially straining the communication resources. 
However, since many of the vehicles will be in free flow, moving at the same speed as the vehicles around them, the measurements that they would communicate are often redundant.
Therefore, it can be beneficial to develop an algorithm that will only activate the potential probe vehicles when their measurements is needed. 
The complete control scheme is depicted in Figure~\ref{fig:blockDiagram}. 


\subsection{Probe vehicle selection}

It is clear that the quality of traffic state reconstruction can only increase if we gain access to more information, i.e. use more probe vehicles.
However, in a situation where the communication channel bandwidth is limited, it might be useful to reconstruct the density with fewer sensors, eliminating the redundant information.
Information about the congestion and stop-and-go waves is particularly important and will significantly improve the control performance, whereas the information about the rest of the road that is in free flow is less impactful.

\begin{figure}
	\centering
	\includegraphics[width=8cm]{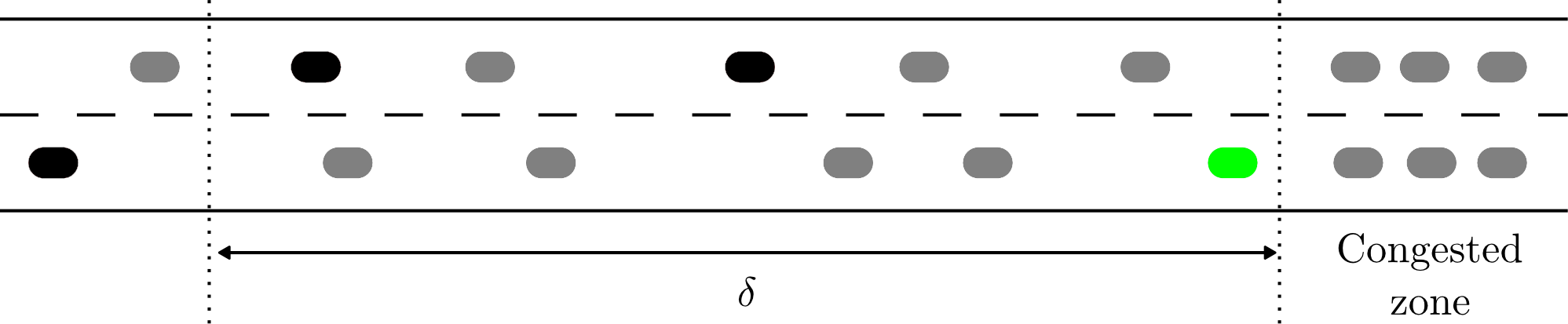}
	\caption{Picture of the road when congestion is detected. It shows the trigger mechanism for activating the dormant CAVs. Grey vehicles are human-driven, without sensing capabilities while the black ones are CAVs that can be switched on to become probe vehicles. Once the active probe vehicle, shown in green, detects the congestion, the black vehicles that are less than $\delta$ away from it activate sensing.}
	\label{fig:selectionScheme}
\end{figure}

There are numerous ways of selecting which CAVs are used as probe vehicles, and the selection will depend on the intended purpose.
Here, we propose a simple adaptive probe vehicle selection scheme outlined in Figure~\ref{fig:selectionScheme}.
Denote by $\mathcal{Y}$ the set of all CAVs. First, we use a subset $\mathcal{Y}^0_s$ of CAVs that are always activated. $\mathcal{Y}_a$ is the set of actuated vehicles ($\mathcal{Y}_a \subset \mathcal{Y}^0_s$) and $\mathcal{Y}_s(t)$ is the set of probe vehicles activated at time $t$.
When the vehicles in $\mathcal{Y}_s^0$ detect any congestion, we activate all vehicles that are less than some predefined distance $\delta$ away from any cell where ${\hat{\rho}_i(t)>\sigma}$, i.e. where we expect to encounter congestion, leading to the following definition:
\begin{equation}
    \mathcal{Y}_s(t) = \mathcal{Y}_s^0 \cup \left\{m \in \mathcal{Y}: \hat{\rho}_i(t)>\sigma, i-i_y^m(t) \in \left[0,\left\lfloor\frac{\delta}{L}\right\rfloor\right]\right\}.
\end{equation}
After these vehicles leave the congestion, they will be deactivated. With this selection scheme, we aim to improve the quality of traffic state reconstruction that is most impactful with regards to the control performance.

In the following section, we will compare this probe vehicle selection scheme to using only the probe vehicles from $\mathcal{Y}^0_s$, as well as with using all CAVs as probe vehicles.


\section{Simulations}
\label{sec:simulations}

\begin{table}
    \centering
    {\normalsize\begin{tabular}{c|c|c}
         Case & Sensors & Actuators \\
         \hline
         No control & $\emptyset$ & $\emptyset$ \\
         Predefined subset of CAVs & $\mathcal{Y}_s^0 \subset \mathcal{Y}$ & $\mathcal{Y}_a \subset \mathcal{Y}_s^0$\\
         Adaptive subset of CAVs & $\mathcal{Y}_s(t) \subset \mathcal{Y}$ & $\mathcal{Y}_a \subset \mathcal{Y}_s^0$\\ 
         All CAVs & $\mathcal{Y}$ & $\mathcal{Y}_a \subset \mathcal{Y}$\\ 
         Full-information & all the road & $\mathcal{Y}_a \subset \mathcal{Y}$ 
    \end{tabular}}
    \caption{Summary of different state reconstruction and control scenarios using CAVs.}
    \label{tab:summary}
\end{table}

\begin{table}
    \centering
    {\normalsize\begin{tabular}{c|c|c}
         $\bar{q}_0 = 3200$ veh/h  & $\sigma=40$ veh/km  & $\alpha=0.25$ \\
         \hline
         $V=100$ km/h  & $W=50$ km/h  & $u_{min}=30$ km/h \\
    \end{tabular}}
    \caption{Simulation parameters.}
    \label{tab:simulationParameter}
    \vspace*{-0.4cm}
\end{table}

\begin{figure*}[t!]
\centering
\begin{subfigure}[t]{0.24\linewidth}
\vspace{-0.05cm}
\includegraphics[height=3.5cm]{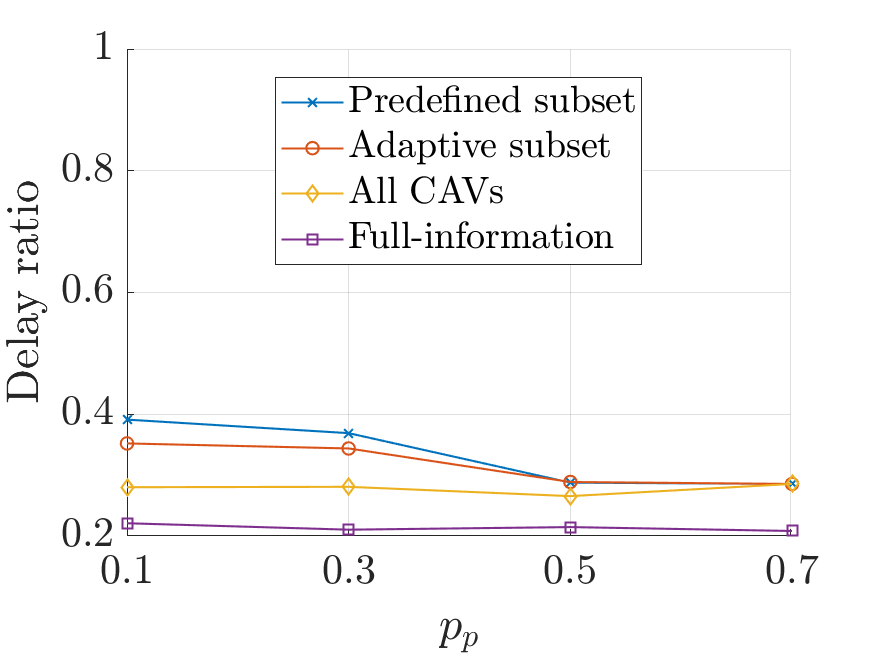}
\caption{$G=0.5$km}
\label{fig:TTT_G100}
\end{subfigure}
\begin{subfigure}[t]{0.23\linewidth}
\vspace{-0.05cm}
\includegraphics[height=3.5cm,trim={0.8cm 0 0 0},clip]{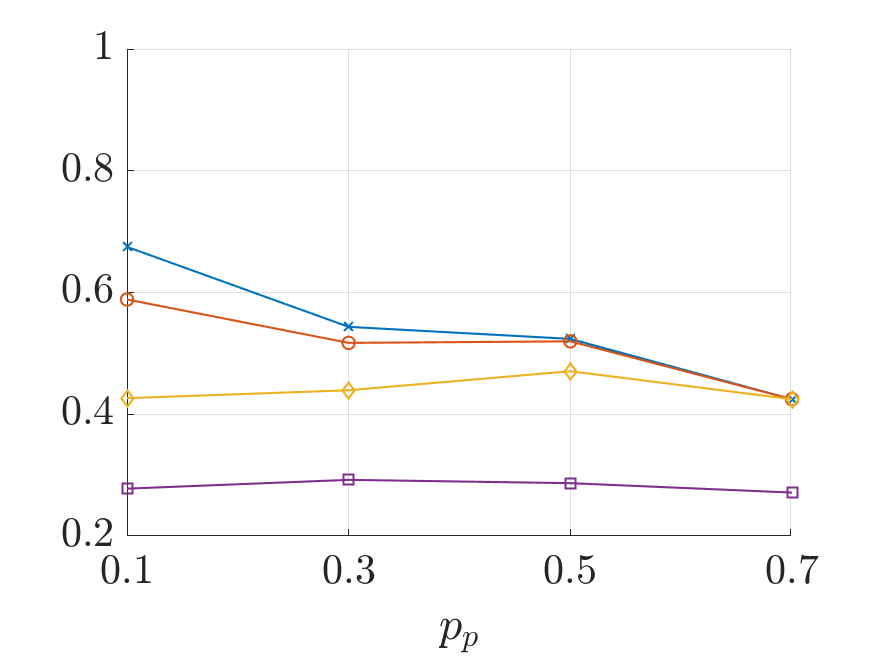}
    \caption{$G=1$km}
    \label{fig:TTT_G200}
\end{subfigure}
\begin{subfigure}[t]{0.23\linewidth}
\vspace{-0.05cm}
\includegraphics[height=3.5cm,trim={0.8cm 0 0 0},clip]{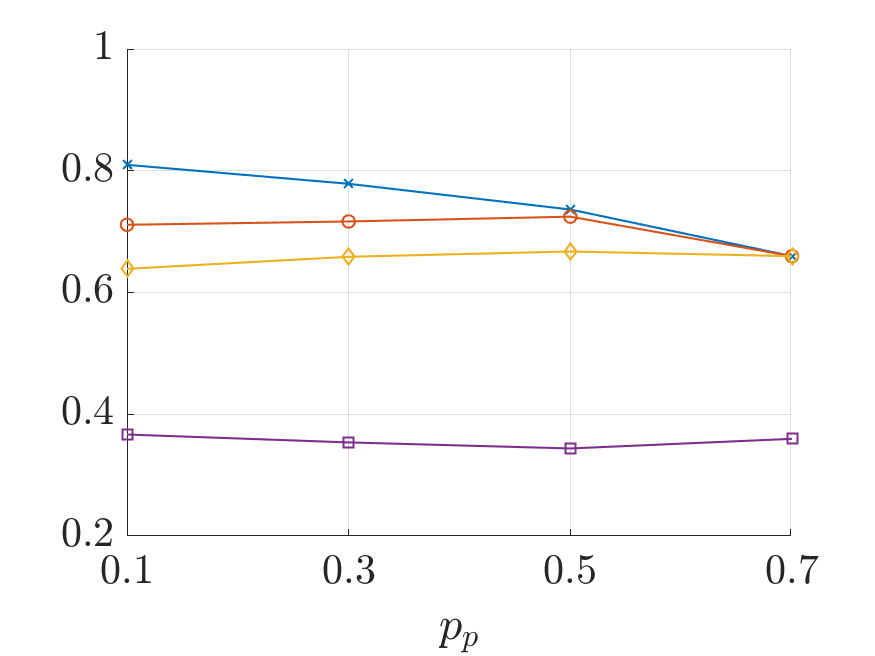}
    \caption{$G=1.5$km}
    \label{fig:TTT_G300}
\end{subfigure}
\begin{subfigure}[t]{0.23\linewidth}
\vspace{-0.05cm}
\includegraphics[height=3.5cm,trim={0.8cm 0 0 0},clip]{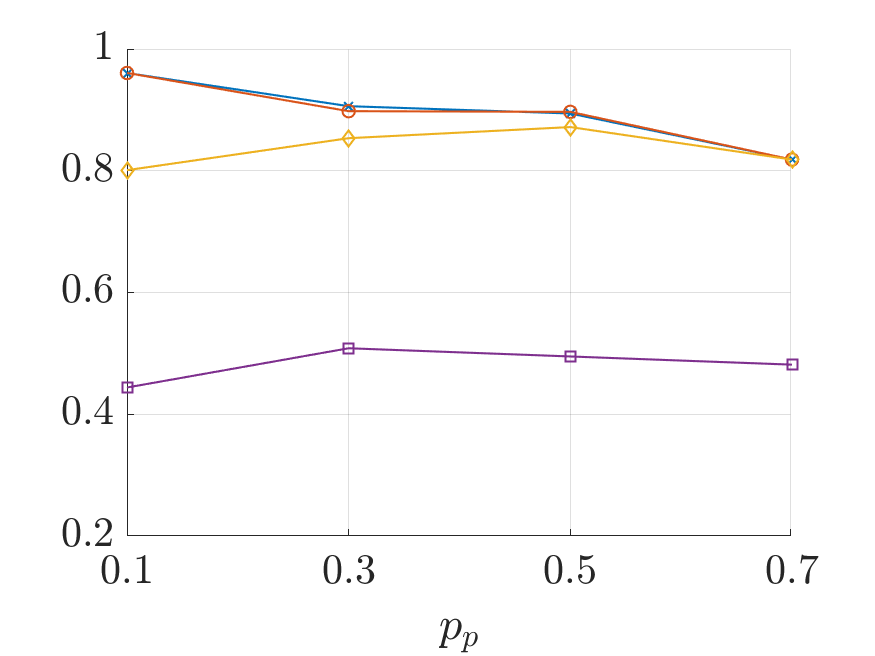}
    \caption{$G=2.5$km}
    \label{fig:TTT_G500}
\end{subfigure}
\caption[]{Median delay compared to the no control case shown for different $p_p$ and $G$ for a given $p_a = 0.3$.}
\label{fig:TTT_G}
\end{figure*}


The efficiency of the proposed traffic density reconstruction scheme and control law is studied in simulations with different parameters.
The initial traffic densities, flow into the road segment and stop-and-go waves entering the road from downstream are randomized and the same values are taken for all parameter combinations.
The arrival of CAVs is modelled as Poisson arrival process with average gap of $G$, and the newly arrived vehicle is a probe vehicle with probability $p_p$, an actuator vehicle with probability $p_a$ (in which case it also acts as a probe vehicle), and an inactive CAV if it is neither a probe nor an actuator vehicle.

We executed $100$ simulations for every different combination of parameter $G$ and $p_s$ values, with ${G \in \left\{0.5\mathrm{km}, 1\mathrm{km}, 1.5 \mathrm{km}, 2.5 \mathrm{km} \right\}}$, ${p_p \in \left\{ 0.1, 0.3, 0.5, 0.7 \right\}}$, with $p_a = 0.3$.
For each simulation run, five control cases were used, as summarized in Table~ \ref{tab:summary}.
The performance metric that we use is the median of the delay ratio,
\begin{equation}
    \frac{TTS - TTS^{\min}}{TTS^{\mathrm{unc}} - TTS^{\min}},
\end{equation}
i.e. the ratio between the increases of $TTS$ in the controlled and uncontrolled case.
The increase is calculated compared to the minimum $TTS^{\min} = \frac{\bar{q}_0}{V}l t_{\mathrm{sim}}$, where ${l = 5}$~km is the length of the simulated road segment and ${t = 1}$~h the duration of each simulation run.

The simulation results are shown in Figure~\ref{fig:TTT_G}.
We can see how increasing $p_p$ affects the control performance for constant $G$ and $p_a$.
Unsurprisingly, we can see that control performance deteriorates as we use less and less information.
When $G=0.5$~km, the full-information control achieves the best performance, eliminating close to $80\%$ of delay, whereas using all CAVs as sensors eliminates around $72\%$ of delay for the same $G$.
In case we are using a subset of CAVs as sensors, the performance will improve as $p_p$ increases, starting from eliminating around $60\%$ and $65\%$ of the delay, using a predefined and adaptive subset of CAVs as sensors respectively, and approaching the performance of the case where we use all CAVs as sensors as $p_p+p_a$ go to $1$, when the same subsets of CAVs are used.
We can also see that the main factor determining the control performance is the average gap between two CAVs.
When $G$ is very low, all control schemes achieve good results, and probe vehicle-based control approaches the full-information control.

\begin{figure*}[t!]
\centering
\begin{subfigure}[t]{0.325\linewidth}
\includegraphics[width=\textwidth]{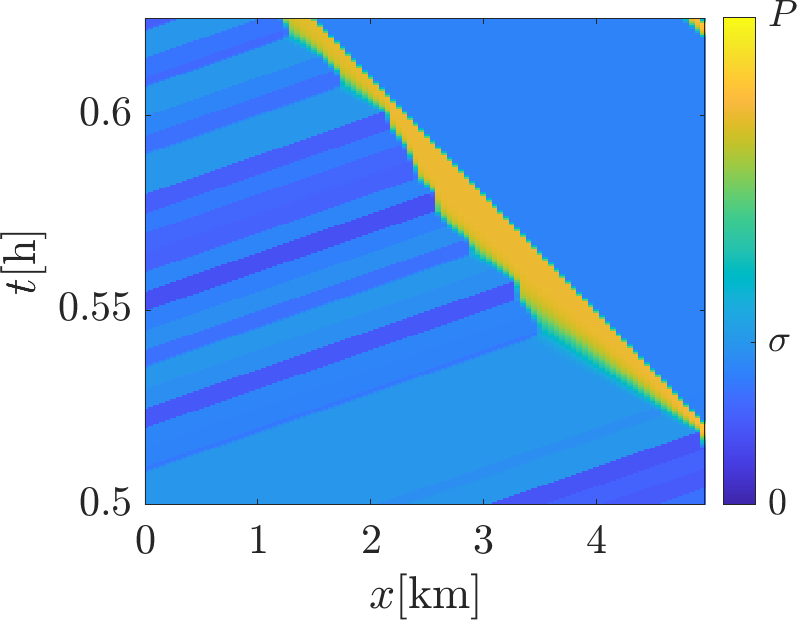}
    \caption{Real density without control.}
    \label{fig:rho_noctrl}
\end{subfigure}
\begin{subfigure}[t]{0.325\linewidth}
\includegraphics[width=\textwidth]{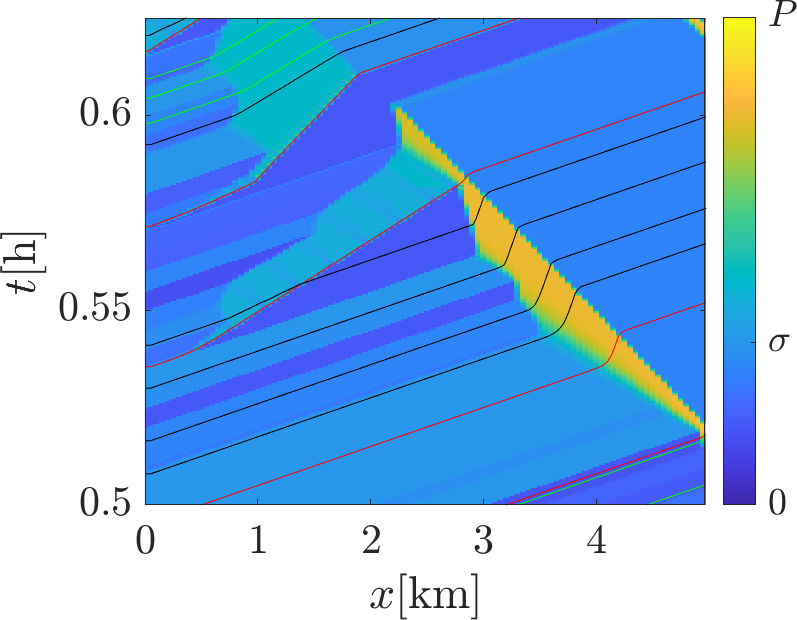}
\caption{Real density with reconstruction-based control using a predefined subset of CAVs.}
\label{fig:rho_predefctrl}
\end{subfigure}
\begin{subfigure}[t]{0.325\linewidth}
\includegraphics[width=\textwidth]{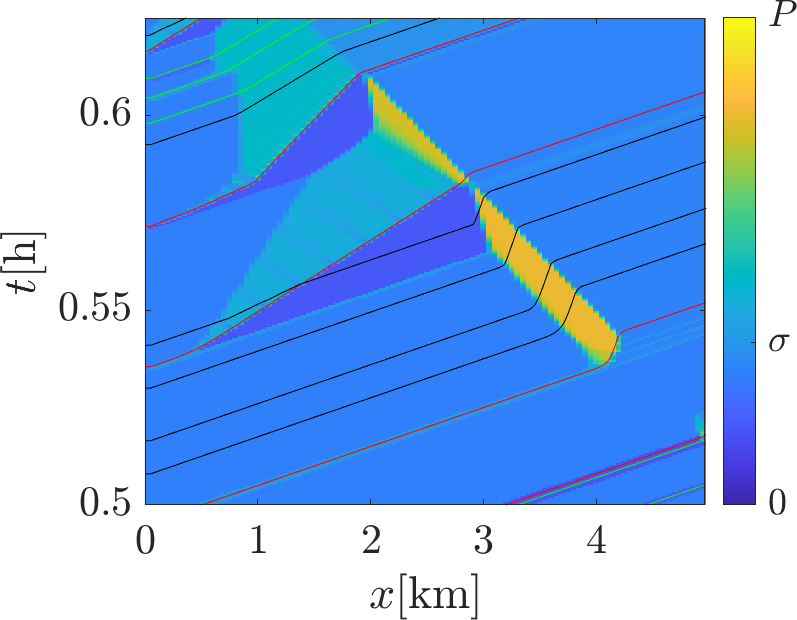}
    \caption{Reconstructed density of situation (b).}
    \label{fig:rhohat_predefctrl}
\end{subfigure}
\\
\vspace{0.2cm}
\begin{subfigure}[t]{0.325\linewidth}
\includegraphics[width=\textwidth]{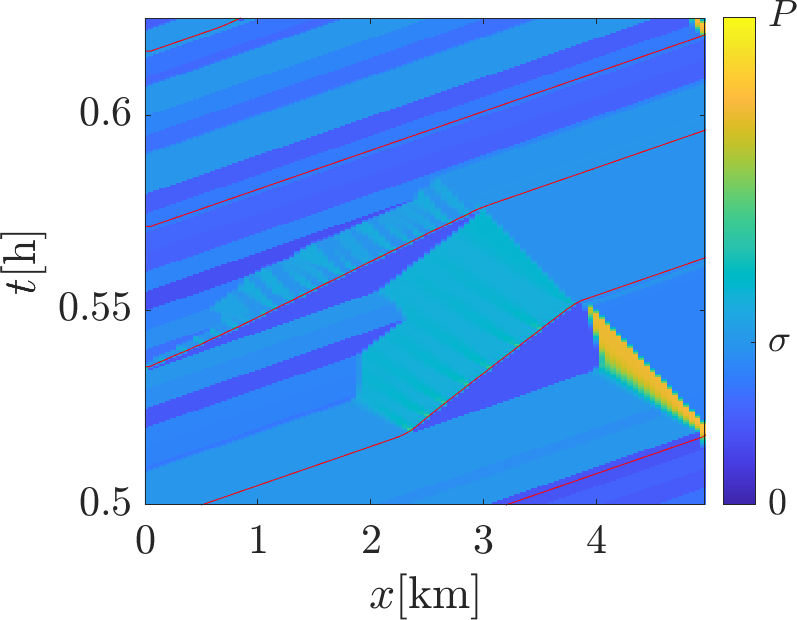}
    \caption{Real density with full information control.}
    \label{fig:rho_fullinfoctrl}
\end{subfigure}
\begin{subfigure}[t]{0.325\linewidth}
\includegraphics[width=\textwidth]{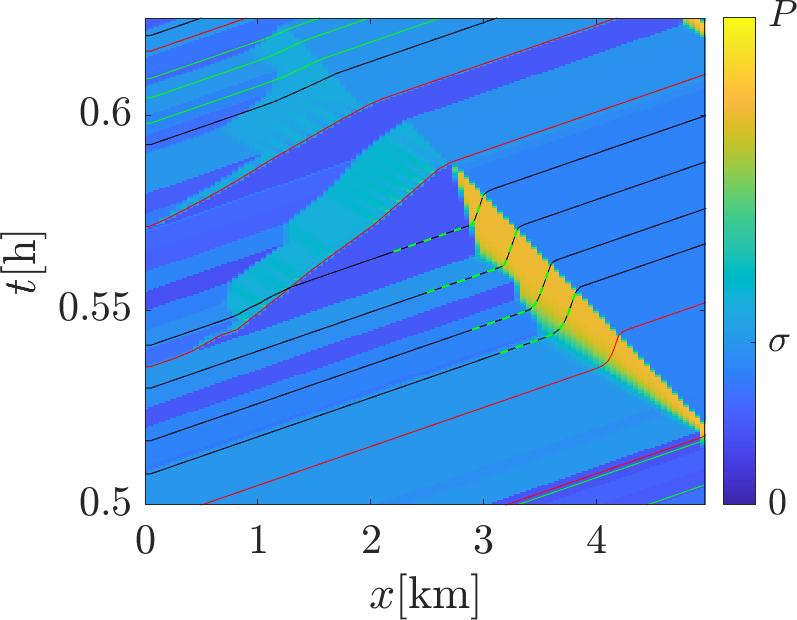}
    \caption{Real density with control using an adaptive subset of CAVs.}
    \label{fig:rho_adaptctrl}
\end{subfigure}
\begin{subfigure}[t]{0.325\linewidth}
\includegraphics[width=\textwidth]{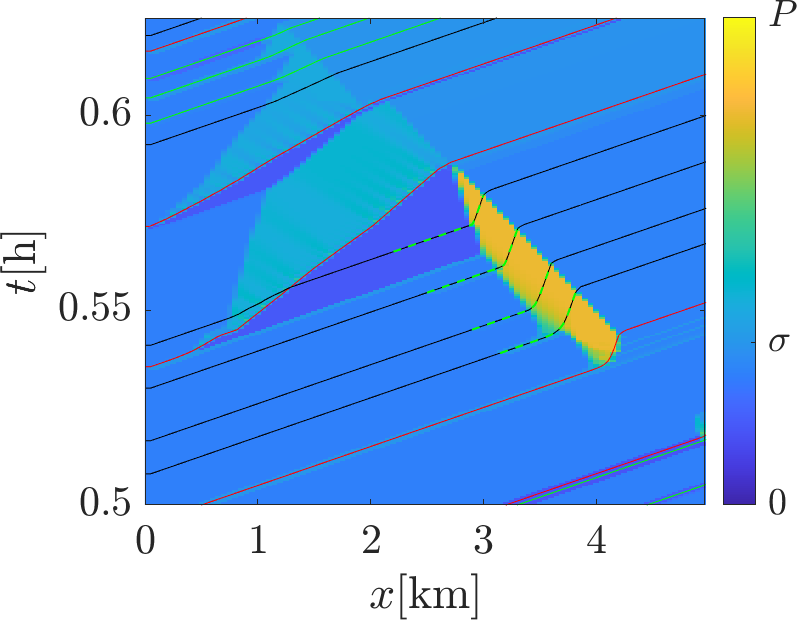}
    \caption{Reconstructed density of situation (e).}
    \label{fig:rhohat_adaptctrl}
\end{subfigure}
\caption[]{Traffic density obtained from a numerical simulation in different cases. The trajectories of inactive CAVs are in black, probe vehicles in green and actuator vehicles in red.}
\label{fig:rho_example}
\end{figure*}

To further illustrate the control and state reconstruction algorithms, in Figure~\ref{fig:rho_example} we show a detail from one of the simulation runs.
The figures show the density profile along the road, with brighter colours representing denser traffic.

The baseline case, with no control, is shown in Figure~\ref{fig:rho_noctrl}, and the full-information control case is shown in Figure~\ref{fig:rho_fullinfoctrl}.
A stop-and-go wave originating from downstream enters the road around ${t=0.52}$~h, and propagates upstream unless dissipated by applying some control action.

Figures~\ref{fig:rho_predefctrl}~and~\ref{fig:rhohat_predefctrl} show the attempt to dissipate the congestion using only the predefined set of CAVs (shown in green and red) as sensors, with the real traffic situation $\rho_i(t)$ shown in Figure~\ref{fig:rho_predefctrl}, and the reconstructed estimation of the traffic state $\hat{\rho}_i(t)$ shown in Figure~\ref{fig:rhohat_predefctrl}.
Around ${t=0.536}$~h, an actuator vehicle runs into the stop-and-go wave, detecting it as it goes through it.
The actuator vehicle upstream reacts by slowing down and restricting the flow.
Four inactive CAVs reach the stop-and-go wave before the actuator vehicle, but since they transmit no information, the control law underestimates the width of the wave and the CAV fails to completely dissipate it.
However, in case we use the proposed adaptive probe vehicle activation, once these four inactive CAVs get close to the congestion, they are temporarily activated, as shown in Figures~\ref{fig:rho_adaptctrl}~and~\ref{fig:rhohat_adaptctrl}  (shown in dashed green).
The additional information corrects the underestimation, and the stop-and-go wave is successfully dissipated.



\section{Conclusion}
\label{sec:conclude}

To conclude, in this paper we explored the effect of using reconstructed traffic state on the efficiency of the control, compared to the full-information case. Both the state reconstruction and actuation are executed using connected automated vehicles. We have seen that the difference between the full-information case and the reconstructed-based control is diminished by introducing more probe vehicles. We also proposed a triggering mechanism which reduces the communication burden, by reducing the number of active probe vehicles, without sacrificing too much control performance.

This preliminary work shows that the Lagrangian state reconstruction should be investigated in more detail since the preliminary results seem very promising. The work conducted here should also be extended to other traffic models. Then, a more refined model could be used with noise in the measurements and external disturbances. Finally, a more formal analysis also needs to be conducted. 

\bibliographystyle{IEEEtran}
\bibliography{biblio}

\end{document}